\documentclass[ aps, pre, superscriptaddress,twocolumn,showpacs]{revtex4}  
\usepackage{graphicx}  
\usepackage{dcolumn}   
\usepackage{bm}        
\usepackage{amssymb}   
\usepackage{scalefnt}

\newcommand{\be}{\begin{equation}}
\newcommand{\ee}{\end{equation}}
\newcommand{\bea}{\begin{eqnarray}}
\newcommand{\eea}{\end{eqnarray}}

\newcommand{\nod}{\noindent}

\newcommand{\ba}{\begin{array}}
\newcommand{\ea}{\end{array}}
\newcommand{\bc}{\begin{center}}
\newcommand{\ec}{\end{center}}

\begin{document}



\title{{\bf Epidemic threshold and control in a dynamic network}}

\author{Michael Taylor}
\affiliation{School of Mathematical and Physical Sciences, Department of Mathematics, University of Sussex, Brighton BN1 9QH, UK}
\author{Timothy J. Taylor }
\affiliation{School of Mathematical and Physical Sciences, Department of Mathematics, University of Sussex, Brighton BN1 9QH, UK}
\affiliation{Centre for Computational Neuroscience and Robotics, University of Sussex, Brighton BN1 9QH, UK}
\author{Istvan Z. Kiss}
\affiliation{School of Mathematical and Physical Sciences, Department of Mathematics, University of Sussex, Brighton BN1 9QH, UK}


\begin{abstract}
In this paper we present  a model describing Susceptible-Infected-Susceptible (SIS) type epidemics spreading on a dynamic contact network with random link activation and deletion where link activation can be locally constrained. We use and adapt a improved effective degree compartmental modelling framework recently proposed by Lindquist et al. [J. Lindquist \emph{et al.}, J. Math Biol. {\bf 62}, 2, 143 (2010)] and Marceau et al. [V. Marceau \emph{et al.},  Phys. Rev. E {\bf 82}, 036116 (2010)]. The resulting set of ordinary differential equations (ODEs) is solved numerically and results are compared to those obtained using individual-based stochastic network simulation. We show that the ODEs display excellent agreement with simulation for the evolution of both the disease and the network, and is able to accurately capture the epidemic threshold for a wide range of parameters. We also present an analytical $R_{0}$ calculation for the dynamic network model and show that depending on the relative timescales of the network evolution and disease transmission two limiting cases are recovered: (i) the static network case when network evolution is slow and (ii) homogeneous random mixing when the network evolution is rapid. We also use our threshold calculation to highlight the dangers of relying on local stability analysis when predicting epidemic outbreaks on evolving networks.
\end{abstract}

\maketitle


\section{Introduction}

The rise in the popularity and relevance of networks as a tool for modelling complex systems is well illustrated by the ever increasing body of research concerned with the spread of diseases within host populations exhibiting non-trivial contact structures \cite{NewmanRev,Albert}. Networks offer an intuitive and relatively simple modelling framework which enables us to relax the strong implicit assumptions of more classical ordinary differential equations (ODE) based approaches and to account for complexities in the contact structure of the host population \cite{NewmanEpiOnNetw,GrossAd,Strogatz,Bansal,EamesKeel02}. This approach has shown that epidemic thresholds not only depend upon the infectiousness of the pathogen, or even simply the mean number of contacts per individual, but also upon  the exact structure of the host population \cite{TrapAnalytic,AndersonSocNet}. In addition to its inherent theoretical value, this paradigm has immediate practical benefits, as the primary role of public health services is to put measures in place to bring diseases below their epidemic threshold. These measures depend heavily upon disrupting the transmission of a disease through vaccination and also more directly through the closure of public services, or even quarantine and curfews in extreme cases. Hence the knowledge of how the structure of the host population is contributing to the spread of a disease would help to increase the efficacy of any intervention \cite{Meyers}.

Despite advances in both rigorous and non-rigorous analysis of networks, a key assumption in many network models is that contacts are fixed for the duration of an epidemic and that the disease propagates with a constant intensity across links. This will not be true for many diseases, especially those with long infectious periods, or diseases that become endemic. Indeed human contact patterns are well described by short repeated events, with individuals having a number of contacts best described by some appropriate time dependent random variable \cite{Read}. Furthermore, individuals and the communities they belong to are likely to change their contact behaviour as a result of natural evolution and endogenous or exogenous perturbations such as a disease outbreak \cite{LilSex}.

 Recently a number of studies have attempted to relax this assumption by allowing the networks to evolve over time by either varying contacts independently of the status of individuals \cite{VolzDyn, VolzThresh} or by explicitly coupling contact activation and deletion to the disease status of individuals \cite{Marceau, GrossAdNet, VanSegAdNet}. Thus, in the latter case, the dynamics of the disease is coupled with the dynamics of the network itself, with both acting as a feedback mechanism for the other \cite{VanSegAdNet, Zanette, GrindEvol}. Many of these studies have built macro ODE-based models that describe the coevolution of networks and the diseases that spread along them \cite{Marceau, GrossAdNet, VanSegAdNet, SchwarzEvolNet}. All these studies confirm that dynamic networks and the coupling between the two dynamics lead to a richer spectrum of behavior than is found for epidemics on static networks.

 A crucial feature of allowing the co-evolution of disease and network is the interplay and feedback between both dynamics, however this interdependence is difficult to measure empirically. The models developed so far mainly use rewiring rules that intuitively make sense given that individuals would have knowledge of the disease states of the rest of the population. However in this paper we move away from these assumptions and we propose a dynamic network model that is based on random link activation-deletion, which would be more relevant for asymptomatic diseases, such as Chlamydia \cite{Low}. Furthermore our dynamic network model is refined by introducing a local constraint on link activation to account for the difference in the magnitude of the number of contacts of a node relative to system size. This dynamic network coupled with the simple Susceptible-Infected-Susceptible (SIS) disease dynamics leads to the full model that will be analysed and discussed. We study this system and explore to what extent a macro ODE-based compartmental model proposed for static networks is flexible enough to be adapted to a dynamic network case. Specifically, we focus on the SIS effective degree model  as described in detail by Lindquist et al.\ \cite{Lind} and also, to our knowledge, proposed by Marceau et al.\ \cite{Marceau} in close succession. Gleeson \cite{Glee} later uses this same modelling framework and demonstrates that the effective degree formulation can be used to model other binary-state dynamics such as Glauber spin dynamics and shows that the ODE model can be used to carry out linear stability type analysis.

Whereas both Lindquist et al.\ \cite{Lind} and Gleeson \cite{Glee} confine themselves to modelling on static contact networks, Marceau et al.\ \cite{Marceau} uses this same improved effective degree formalism to explore SIS disease dynamics on adaptive networks. In this model the number of links in the network is fixed but the susceptible individuals can replace links to infectious neighbours with links to other randomly chosen susceptible individuals, as originally proposed by Gross et al.\ \cite{GrossAdNet}. Our proposed model also uses SIS type epidemics on dynamic networks, but unlike Marceau our model allows for the random activation and deletion of links over time. As such not only the network topology will evolve evolve and change over time, but also the number of links. This modified dynamic effective degree model is also governed by a closed set of ODEs, which is then solved and compared to results from individual based simulations and its ability to accurately predict the epidemic threshold over a range of parameters is investigated. We also derive an analytical $R_{0}$ calculation that describes the stability of the disease-free equilibrium and we discuss the limitations of such a calculation in the light of having a dynamically active and evolving contact network.

\begin{figure}[!t]
\includegraphics[scale=1.0]{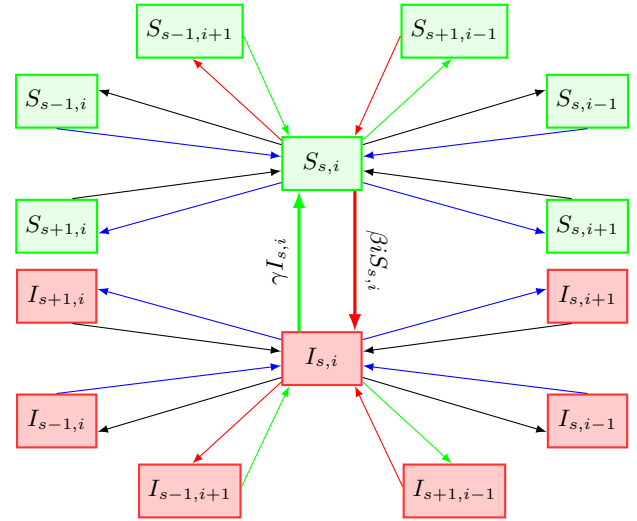}

\caption{\label{chart}
(Colour online) Flow chart showing transitions in the dynamic SIS effective degree model. The directed red (gray), green (light gray), blue (dark gray) and black lines represent changes in state of an individual via infection, recovery, link creation and link deletion respectively. The thick lines represent changes to the individual, and thin lines represent changes to that individual's immediate neighbourhood. In relation to nodes of type $X_{si}, X \in \{S,I\}$, infection of neighbours occurs at rate $sG_{X}$, recovery of neighbours at rate $\gamma i $, creation of a susceptible (infectious) link at rate $\alpha(M-(s+i))P_{S(I)}$ and deletion of a susceptible (infectious) link at rate $\omega s (i)$, where:\\
$G_{S}=\beta \frac{\sum_{k=1}^{M}\sum_{j+l=k}jlS_{jl}}{\sum_{k=1}^{M}\sum_{j+l=k}jS_{jl}}$, $G_{I}=\beta \frac{\sum_{k=1}^{M}\sum_{j+l=k}l^{2}S_{jl}}{\sum_{k=1}^{M}\sum_{j+l=k}jI_{jl}}$\\ and $P_{X}=\frac{\sum_{k=0}^{M}\sum_{j+l=k}(M-(j+l))X_{jl}}{\sum_{k=0}^{M}\sum_{j+l=k}(M-(j+l))(S_{jl}+I_{jl})}$.}
\end{figure}


\section{The model}

 Linquist et al. \cite{Lind} and Marceau et al. \cite{Marceau} use different notation to describe the same modelling framework. For consistency, in this paper we follow the notation used by the former throughout. The effective degree modelling approach for SIS type disease dynamics \cite{Lind} not only categorizes the disease state of each individual as susceptible ($S$) or infected ($I$) but also describes the state of their immediate neighbourhood. This is achieved by keeping track of the number of susceptible and infected neighbours that belongs to a given node. For example, $S_{si}$ represents the number of susceptible individuals that have $s$ susceptible and $i$ infected neighbours. This gives rise to more states and equations than would be seen in a standard pairwise model, where equations are given at the population level for all types of singles and pairs \cite{KeeEa}. For example if a $S_{si}$ type node became infected via one of its $i$ infectious neighbours, this individual would move to state $I_{si}$ as only the status of the node itself is changing. However, if one of the $i$ infected neighbours of an $S_{si}$ type node recovered then the node would enter the $S_{s+1,i-1}$ class, whereas infection of one of the $s$ neighbouring susceptible nodes moves the $S_{si}$ type node into the $S_{s-1,i+1}$ class.

Lindquist et al. \cite{Lind} defined $\gamma$ to be the per node recovery rate, $\beta$ the per link infection rate and $M$  the maximum nodal degree of a network with $N$ nodes. They then derived the following system of $\sum_{k=1}^{M}2(k+1)=M(M+3)$ equations:

\bea
\dot{S_{si}}&=&-\beta i S_{si} + \gamma I_{si} +\gamma [(i+1)S_{s-1,i+1}-iS_{si}]\label{statSsi} \\
&&+ \beta \frac{\sum_{k=1}^{M}\sum_{j+l=k}jlS_{jl}}{\sum_{k=1}^{M}\sum_{j+l=k}jS_{jl}}[(s+1)S_{s+1,i-1} -sS_{si}], \nonumber \\
\dot{I_{si}}&=&\beta i S_{si} - \gamma I_{si} +\gamma [(i+1)I_{s-1,i+1}-iI_{si}] ]\label{statIsi}\\
&&+ \beta \frac{\sum_{k=1}^{M}\sum_{j+l=k}l^{2}S_{jl}}{\sum_{k=1}^{M}\sum_{j+l=k}jI_{jl}}[(s+1)I_{s+1,i-1} -sI_{si}], \nonumber
\eea

\nod for $\{(s,i):s, i\ge 0, 1 \le k=s+i \le M \}$. This is the SIS effective degree model for a \emph{static} contact network.

In oder to adapt this model to describe SIS dynamics on a \emph{dynamic} contact network, we introduce two new parameters: $\omega$, the per link deletion rate and $\alpha$, the per non-link, or more precisely the per \emph{potential} link creation rate. These rates could also be made to be link-type dependent, i.e. $\omega_{SI}$ would be the per SI link deletion rate. For the dynamic network cae, the system size will increase slightly from $M(M+3)$ to $\sum_{k=0}^{M}2(k+1)=(M+1)(M+2)$ equations to account for nodes of the type $X_{0,0}$ where $X \in \{S,I\}$. In the static case, these nodes were dynamically unimportant as they could neither infect nor become infected by other nodes. However in the dynamic model, they could connect to other nodes in the system and so enter states $X_{1,0}$ or $X_{0,1}$ depending on the state of the node with which they have just formed a new link.

The total number of links in the system at time $t$, $\Lambda(t)$, and potential links, $\Phi(t)$ can easily be calculated from the effective degree formulation as

\bea
\Lambda(t)&=&\sum_{k=0}^{M}\sum_{j+l=k}(j+l)(S_{jl}+I_{jl}), \nonumber \\
\Phi(t)&=&\sum_{k=0}^{M}\sum_{j+l=k}\left(M-(j+l)\right)(S_{jl}+I_{jl}) ], \nonumber
\eea

%
\nod with the mean nodal degree given by $\langle k(t)\rangle =\frac{\Lambda(t)}{N}$. At the equilibrium, $\alpha \Phi = \omega \Lambda$ which gives us the mean nodal degree:

\be
\langle k \rangle^{*}=\frac{\alpha}{\alpha + \omega}M. \label{degEq}
\ee

\nod Note that Eq.\ (\ref{degEq}) does not depend on the system size, $N$, but rather on the maximum nodal degree, $M$. This is important because in the static model, $M$ is simply given by the node or nodes with the highest degree whilst in the dynamic case, however, $M$ can be considered as a carrying capacity, whereby no node can have more than $M$ links. This subtle but important difference means that in the dynamic case, $M$ itself can be regarded as a parameter which controls the potential level of network saturation.

\begin{figure}[!t]
\includegraphics[scale=0.7]{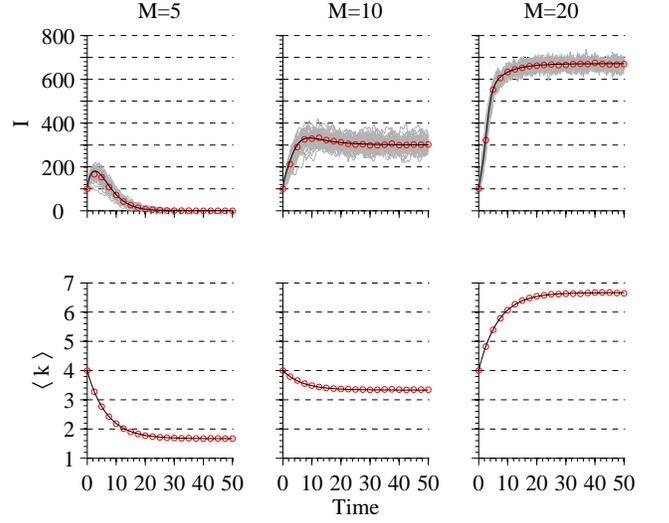}
\caption{\label{fig1} (Colour online) Time evolution of $I(t)=\sum_{k=0}^{M}\sum_{j+l=k}I_{jl}(t)$ and $\langle k \rangle(t) =\frac{\Lambda(t)}{N}$ for three different values of $M$. Results from the ODE are given by solid lines and those from simulation by points. In all cases $N=1000$, $I_{0}=100$, $\alpha=0.05$, $\omega=0.1$, $\beta=0.5$ and $\gamma=1$. The initial network is a regular random graph with $k=4$. In each case, mean values from the stochastic simulations were found by averaging over 100 repetitions, with the individual realisations plotted in grey.}
\end{figure}

When adding the terms that govern link creation and deletion to Eqs.\ (\ref{statSsi}) and (\ref{statIsi}) it is far simpler to construct the terms that govern deletion of existing links than those for the creation of new links. Links to nodes of type $X_{si}$ where $X \in \{S,I\}$ are cut at a rate proportional to their degree, so individuals will leave $X_{si}$ through link deletion at a rate $\omega(s+i)$ and will either enter the $X_{s-1,i}$ or $X_{s,i-1}$ classes depending on the state of the nodes to which they were previously connected. Similarly individuals can enter state $X_{si}$ if they were in states $X_{s,i+1}$ or $X_{s+1,i}$ and a link to an infected or susceptible node was deleted respectively.

When creating new links to nodes of type $X_{si}$, there are $M-(s+i)$ stubs remaining, so nodes will transition out of this state at a rate $\alpha (M-(s+i))$ and will either enter the $X_{s+1,i}$ or $X_{s,i+1}$ classes depending on the state of the node to which they have just connected. The rate at which nodes enter the $X_{si}$ class from either $X_{s-1,i}$ or $X_{s,i-1}$ depends not only on the number of stubs still available in the node in question, but also on the probability that the newly created link attaches to a node of state $S$ or $I$ respectively. So nodes enter $X_{si}$ from $X_{s-1,i}$ at the rate $\alpha P_{S} (M-(s-1+i))$, and nodes enter $X_{si}$ from $X_{s,i-1}$ at rate $\alpha P_{I} (M-(s+i-1))$, where $P_{X}=\frac{\sum_{k=0}^{M}\sum_{j+l=k}(M-(j+l))X_{jl}}{\sum_{k=0}^{M}\sum_{j+l=k}(M-(j+l))(S_{jl}+I_{jl})}, X \in \{S,I\}$ is the probability of picking an available stub belonging to nodes of type $X$ where $X \in \{S,I\}$. The full set of transitions captured by this model is shown in Fig. \ref{chart}.

The addition of these terms to Eqs.\ (\ref{statSsi}) and (\ref{statIsi}) transforms the SIS effective degree model for a static network into one that captures the spread of SIS type diseases on a dynamic contact network and is described by the following system of $(M+1)(M+2)$ equations:

\bea
\dot{S_{si}}&=&-\beta i S_{si} + \gamma I_{si} +\gamma [(i+1)S_{s-1,i+1}-iS_{si}] \label{dynSsi} \\
&&+ \beta \frac{\sum_{k=0}^{M}\sum_{j+l=k}jlS_{jl}}{\sum_{k=0}^{M}\sum_{j+l=k}jS_{jl}}[(s+1)S_{s+1,i-1} -sS_{si}] \nonumber \\
&&-\omega[(s+i)S_{si} -(i+1)S_{s,i+1}-(s+1)S_{s+1,i}] \nonumber \\&&-\alpha(M-(s+i))S_{si}   +\alpha(M-(s-1+i))P_{S}S_{s-1,i}, \nonumber \\
&& +\alpha(M-(s+i-1))P_{I}S_{s,i-1} \nonumber \\
\dot{I_{si}}&=&\beta i S_{si} - \gamma I_{si} +\gamma [(i+1)I_{s-1,i+1}-iI_{si}] \label{dynIsi} \\
&&+ \beta \frac{\sum_{k=1}^{M}\sum_{j+l=k}l^{2}S_{jl}}{\sum_{k=1}^{M}\sum_{j+l=k}jI_{jl}}[(s+1)I_{s+1,i-1} -sI_{si}] \nonumber \\
&&-\omega[(s+i)I_{si} -(i+1)I_{s,i+1}-(s+1)I_{s+1,i}] \nonumber \\&&-\alpha(M-(s+i))I_{si}
 +\alpha(M-(s-1+i))P_{S}I_{s-1,i} \nonumber \\
&& +\alpha(M-(s+i-1))P_{I}I_{s,i-1}, \nonumber
\eea

\nod for $\{(s,i):s,i \ge 0, 0 \le k=s+i \le M \}$. This system is the \emph{dynamic} SIS effective degree model.

\section{Calculating the disease threshold}

For the static case, Lindquist et al. \cite{Lind}  used the next generation matrix approach \cite{Diek} to calculate the disease threshold to be

\begin{equation}
{\cal R}_{0}=\rho(FV^{-1})=\frac{\beta}{\sum_{k=1}^{M}kS_{k,0}}\sum_{k=1}^{M}v_{k}^{T}V_{k}^{-1}u_{k}.\label{R0stat}
\end{equation}

In this approach, Eqs. (\ref{dynSsi}) and (\ref{dynIsi}) are linearized at the disease-free equilibrium (DFE) and the Jacobian at the DFE is written as $F-V$. In this formulation, F accounts for transitions from disease-free states to disease states (in the static case only the transition from $S_{s,0}$ to $S_{s-1,1}$ needs to be considered) and $V$ accounts for transitions between different disease states. The spectral radius, $\rho$, the leading eigenvalue of $FV^{-1}$, gives $R_{0}$ and describes the stability of the DFE. If $R_{0}<1$ the DFE is stable and no epidemic will occur, but if $R_{0}>1$ the DFE is unstable and the infectious agent can spread through the population.

We can calculate $F$ in the dynamic case by noting that the same $S_{s,0}$ to $S_{s-1,1}$ type transitions can still occur, but in addition nodes can enter the disease states by linking to an infected node, namely $S_{s,0}$ to $S_{s,1}$ transitions. If we introduce a subscript $s$ to denote the static version of the next generation matrix, so the static version of $F$ is called $F_{s}$ and so on, we have

\begin{figure}[!t]
\includegraphics[scale=0.7]{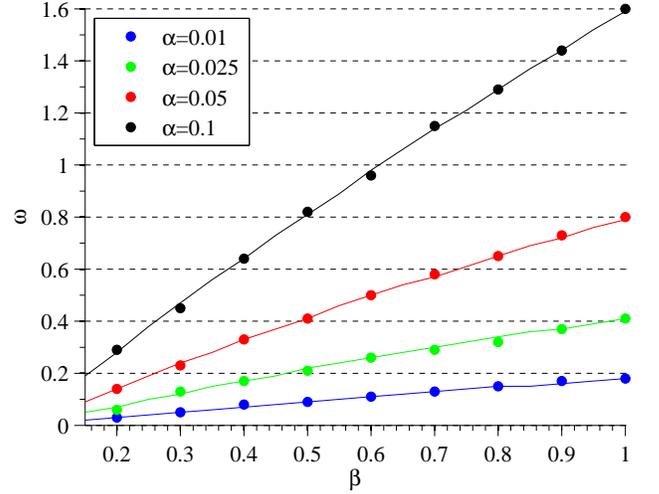}
\caption{\label{fig2} (Colour online) Epidemic threshold plot in the $(\beta, \omega)$ parameter space for four distinct values of $\alpha$.  Results from the ODE are given by solid lines and those from simulation by solid points. In each case, $N=1000$, $I_{0}=10$, $M=20$ and $\gamma=1$. The initial network is a regular random graph with $k=4$.}
\end{figure}

\begin{equation}
F_{s}=\frac{\beta}{\sum_{k=0}^{M}kS_{k,0}} \left[ \begin{array}{c}u_{s_{0}}\\u_{s_{1}}\\ \vdots \\ u_{s_{M}} \end{array} \right] \left[ \begin{array}{cccc} v_{s_{0}}^{T}&v_{s_{1}}^{T}&\dots&v_{s_{M}}^{T} \end{array} \right],
\end{equation}

\nod where $u_{s_{k}}$ and $v_{s_{k}}$ are $(2k+1)\mbox{ x }1$ vectors. The $u_{s_{k}}$ vectors have $kS_{k,0}$ as their first entry and zeros elsewhere and the $v_{s_{k}}$ vectors have their first $(k-1)$ entries equal to $(k-1),2(k-2),\dots,s(k-s),\dots,(k-1)$ and zeros elsewhere. This is almost identical to the $F$ matrix constructed by Lindquist et al., but is augmented by $u_{s_{0}}$ and $v_{s_{0}}$ to account for the new disease state, $I_{0,0}$, and the summation starts at $k=0$ rather than $k=1$,

We now introduce a new subscript $d$ to describe the new transitions that are only possible in the dynamic model. Hence a new $F$ matrix, $F_{d}$, is created, which has exactly the same dimensions as $F_{s}$, and is given by

\begin{equation}
F_{d}=\frac{\alpha}{\sum_{k=0}^{M}(M-k)S_{k,0}} \left[ \begin{array}{c}u_{d_{0}}\\u_{d_{1}}\\ \vdots \\ u_{d_{M}} \end{array} \right] \left[ \begin{array}{cccc} v_{d_{0}}^{T}&v_{d_{1}}^{T}&\dots&v_{d_{M}}^{T} \end{array} \right].
\end{equation}.

Here, $u_{d_{k}}$ is again a $(2k+1)\mbox{ x }1$ vector with the first entry equal to $(M-(k-1))S_{k-1,0}$ and all other entries equal to zero. In the case where $k=0$, $u_{d_{0}}=(0)$. In addition, $v_{d_{k}}$ is the same size as $u_{d_{k}}$ and the first $k$ entries are equal to zero, with the remaining $k+1$ entries equal to $M-k$. The final $F$ matrix that captures all the possible transitions in the dynamic effective degree model is found by taking a linear sum of the two, namely $F=F_{s}+F_{d}$.

\begin{figure}[!t]
\includegraphics[scale=0.65]{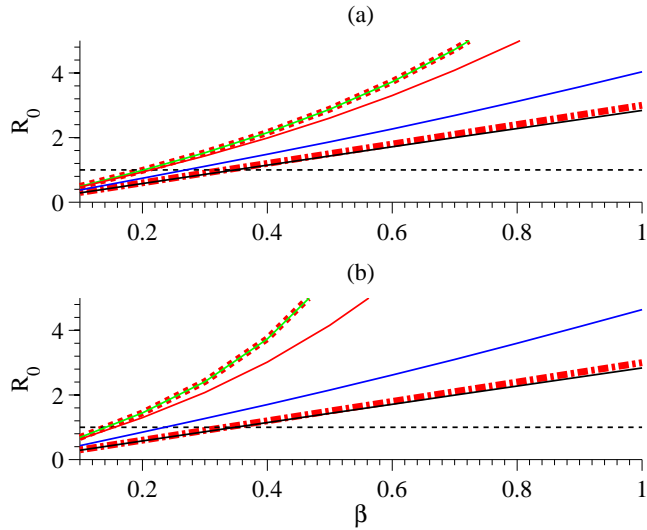}
\caption{\label{fig3} (Colour online) Threshold stability in the $(\beta,R_{0})$ space with $\gamma=1$, $M=20$ and $\langle k \rangle^{*}=3$ for (thin solid lines, in order from top to bottom) $\alpha=10^{-4}$ (green), $\alpha=10^{-2}$ (red), $\alpha=10^{-1}$ (blue) and $\alpha=10$ (black). In (a) the initial network is a regular random graph with $k=6$ and in (b) the initial degree distribution is negative binomial with $\langle k \rangle=6$ and $\sigma^{2}=12$. In each case, $\omega=\alpha\frac{M-\langle k \rangle^{*}}{\langle k \rangle^{*}}$. The thick short-dashed red line is the theoretical value of $R_{0}$ for a static network, and the thick red dash-dotted line is the mean field limit $R_{0}=\frac{\beta}{\gamma}\langle k \rangle^{*}$.}
\end{figure}

As with the static case, the $V$ matrix is constructed through careful book-keeping, which can be done through iterative routines. In the static case, as the nodes have fixed degree, $V_{s}$ is a block diagonal matrix with $V_{s}=V_{s_{1}}\oplus V_{s_{2}}\oplus...\oplus V_{s_{M}}$. For the dynamic model, $V_{d}$ will be a block tri-diagonal matrix, as state transitions can now also occur by nodes gaining or loosing a link. In addition, the extra disease state $I_{0,0}$ needs to be considered,  and $V$ will now also depend upon $\alpha$ and $\omega$ as well as $\beta$ and $\gamma$. Once $F=F_{s}+F_{d}$ and $V=V_{d}$ are constructed, the leading eigenvalue or $R_{0}$ is computed numerically.

\section{Results and discussion}

As shown in Fig.\ \ref{fig1}, the ODEs given by Eqs.\ (\ref{dynSsi}) and  (\ref{dynIsi}) closely capture the time evolution of an epidemic as predicted by stochastic simulations. The only parameter that is varied in Fig.\ \ref{fig1} is $M$, and it is interesting to note the effect it has on the evolution of the disease. As per Eq.\ (\ref{degEq}), the mean nodal degree at equilibrium is dependent on $M$, and hence, given the same initial network configuration and values of $\alpha$ and $\omega$, the network either loses or gains links as the system evolves. Thus varying the carrying capacity alone leads to different outcomes depending on whether the network can reach a level of connectedness that allows an epidemic to spread and become established. Allowing $M$ to become an active model parameter that is able to control the outcome of an epidemic has potentially interesting real world implications. The number of contacts per person is a natural, countable property unlike the other model parameters, such as $\omega$, which are more difficult to infer. Therefore local constraints that limit the maximum number of contacts per person could be potentially used as a metric when promoting safe behaviour at a population level in the event of an outbreak or other public health crisis.

\begin{figure}[!t]
\includegraphics[scale=0.65]{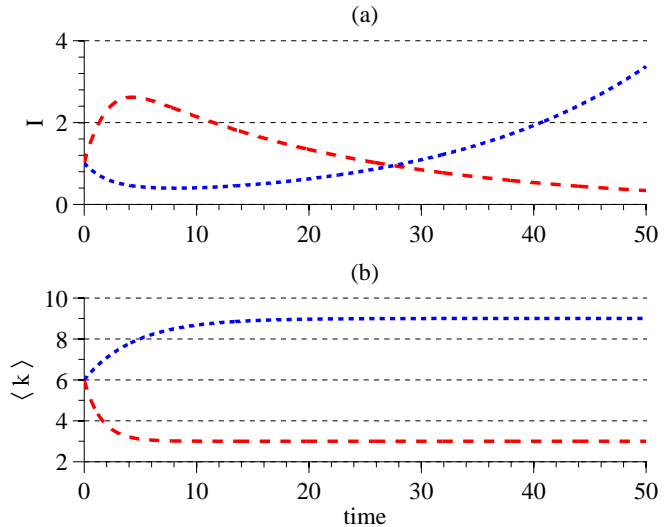}
\caption{\label{fig4} (Colour online) Time evolution of $I$ and $\langle k \rangle$ with $\gamma=1$, $M=20$, $\alpha=0.1$, $\omega=\alpha\frac{M-\langle k \rangle^{*}}{\langle k \rangle^{*}}$ and an initial regular random network with $k=6$. The two cases illustrated above correspond to: $\langle k \rangle^{*}=3$ and $\beta=0.35$, giving $R_{0}\approx 1.29$ (red long-dashed line) and $\langle k \rangle^{*}=9$ and $\beta=0.125$, giving $R_{0}\approx 0.77$ (blue short-dashed line).}
\end{figure}

In Fig.\ \ref{fig2}, for a given value of $\alpha$, $M$ and $\beta$, the epidemic threshold has been calculated from the ODEs in terms of $\omega$ and compared to that predicted by simulations. The agreement is excellent and this is strong evidence that the dynamic effective degree model accurately captures the evolution of an epidemic on a network with random link creation and deletion. When considering the $(\beta,\omega)$ parameter space used for the threshold plot in Fig.\ \ref{fig2}, there are three distinct regions that are worth noting. Firstly, given an initial starting network, it is possible to calculate the threshold value of $\beta$ in the static network case. For the regular random graph with $k=4$ used here, that value is $\beta^{*}\approx 0.36$. For values of $\beta<0.36$, the relative time scales of disease and network evolution are crucial in determining whether or not an epidemic will occur. In this situation, the network needs to quickly evolve to become more densely connected in order for there to be an outbreak. The second area of interest is when the disease is highly infectious and as a result requires a high value of $\omega$  to drive the epidemic below threshold.  Indeed, if the disease parameters $\beta$ and $\gamma$ are fixed then the only way of affecting the outcome of an epidemic is through changing the network structure, i.e. reducing the number of links or the variance. Hence, for a fixed $\alpha$ and $M$, a value of $\beta$ can be chosen large enough so that the minimum value of $\omega$ needed to reduce the connectivity of the network sufficiently to stop an outbreak (see Fig.\ref{fig2}), gives $\langle k \rangle^{*}<2$ as can be calculated from Eq.\ \ref{degEq}. If a network has $\langle k \rangle^{*}<2$ then it becomes fragmented, with many nodes becoming unconnected. In these situations, the value $\omega$ needed to prevent an epidemic virtually destroys the network. In terms of real world implications, a large value of $\omega$ could correspond to a situation of strict quarantine and curfew whereby links between individuals are kept to a minimum. In between these two cases lies a region within which an epidemic would take hold naturally, given the initial network, but which can be prevented by a value of $\omega$ that leaves the network well connected.

 In Fig.\ \ref{fig3}, we show analytical values of $R_{0}$ for a range of values of $\beta$ and $\alpha$. It is worth noting that two limiting cases are recovered when the timescale of the network dynamics is fast and slow relative to the timescale of the disease dynamics. The thick short-dashed red line shows $R_{0}$ calculated for a static network, as proposed by Lindquist et al. \cite{Lind} and given in Eq.\ \ref{R0stat}, and this is exactly followed by results from our dynamic $R_{0}$ calculation when the network dynamics are set to be much slower than the disease dynamics. The other extreme is shown by the thick dash-dotted red line, and is the value of $R_{0}$ that results form the classic mean-field calculation $R_{0}=\frac{\langle k \rangle \beta}{\gamma}$. The time evolution of $\langle k \rangle$ is given by $\dot{\langle k \rangle}=\alpha(M-\langle k \rangle)-\omega \langle k \rangle$ but, when the network dynamics is fast, the equilibrium network distribution, and hence $\langle k \rangle^{*}$, is approached much quicker than the epidemic timescale and hence a value of $\langle k \rangle = \langle k \rangle^{*}$ as given by Eq.\ \ref{degEq} can be used. This limit is closely matched by results from our dynamic $R_{0}$ calculation when the network dynamics are rapid compared to disease transmission as shown in Fig.\ \ref{fig3}.

Although Fig.\ \ref{fig3} demonstrates the accuracy of our analytical $R_{0}$ calculation,  Fig.\ \ref{fig4} highlights two example cases where the long term epidemic outcomes are the opposite of what is predicted by $R_{0}$. In the cases $R_{0}<1$ (blue short-dashed curve) and $R_{0}>1$ (red long-dashed curve) the system settles to an endemic and to a disease free equilibrium respectively, due to the different ways the networks evolve. Given that $R_{0}$ is based on a local stability analysis, it can only incorporate the immediate next-generation effects of random link activation and deletion, and cannot account for long term changes to the network structure. It is well established in the literature (see, for example Li et al. \cite{JingLi}) that $R_{0}$ is of limited value when used as a predictor, and even for static networks needs to be used with care. Our results add weight to this argument, and we show that when dealing with disease spreading through dynamic contact networks the use of $R_{0}$ as any kind of predictor on long term disease evolution should be met with some degree of caution.

In summary, this paper has proposed an effective degree model for epidemics on dynamic networks with random link activation and deletion, where activation is locally constrained. We have shown that this  model agrees extremely well with results obtained from stochastic simulations, and as such can reliably be used for the analytical and semi-analytical study of coupled disease and network dynamics. We have shown how a local constraint limiting the number of contacts per individual can be used to control and prevent the outbreak of an epidemic in this dynamic model. We have also proposed an analytical calculation of $R_{0}$, but also demonstrated the limited value of threshold stability analysis in predicting the evolution of a disease in a dynamic contact network. In future work, this  model can be adapted and extended to account for individuals cutting and creating links with knowledge of the state of others in the population, i.e., link-type dependent network dynamics. This two-way feedback will lead to more sophisticated network properties such as degree correlations, high clustering or even network fragmentation. In such cases ODE models need to be used with care, making sure that the agreement with simulations remains valid. Besides modelling epidemics, this framework could also be used to study the spread of information, beliefs and new ideas within populations, and as such could have implications across a wide range of disciplines beyond the mathematical biology and physics communities.

\noindent {\bf Acknowledgments}
Michael Taylor acknowledges support from EPSRC (DTA grant).
Timothy J. Taylor acknowledges support from MRC.


\end{document}